\documentclass[number,12pt]{elsarticle}   % Elsevier EABE template
\usepackage{graphicx}
\usepackage{natbib}
\usepackage{amssymb}
\usepackage{amsmath}
\usepackage{subfig}
\usepackage{ifthen}
\usepackage{setspace}
\usepackage{tabularx,ragged2e,booktabs,caption}

\providecommand\bnabla{\boldsymbol{\nabla}}
\providecommand\bcdot{\boldsymbol{\cdot}}
\newcommand\p{\ensuremath{\partial}}

\begin{document}

\begin{frontmatter}



\title{A robust and non-singular formulation of the boundary integral method for the potential problem}


\author[1]{Qiang Sun}
\author[2]{Evert Klaseboer}
\author[1]{Boo Cheong Khoo}
\author[1,2,3,4]{Derek Y. C. Chan\corref{cor1}}
\cortext[cor1]{Email: d.chan@unimelb.edu.au}

\address[1]{Department of Mechanical Engineering, National University of Singapore, 10 Kent Ridge Crescent, 119260, Singapore}
\address[2]{Institute of High Performance Computing, 1 Fusionopolis Way, 138632, Singapore}
\address[3]{Department of Mathematics and Statistics, The University of Melbourne, Parkville 3010 VIC Australia}
\address[4]{Department of Chemistry and Biotechnology, Swinburne University of Technology, Hawthorn 3122 VIC Australia}

\begin{abstract}
A non-singular formulation of the boundary integral method (BIM) is presented for the Laplace equation whereby the well-known singularities that arise from the fundamental solution are eliminated analytically. A key advantage of this approach is that numerical errors that arise due to the proximity of nodes located on osculating boundaries are suppressed. This is particularly relevant in multi-scale problems where high accuracy is required without undue increase in computational cost when the spacing between boundaries become much smaller than their characteristic dimensions. The elimination of the singularities means that standard quadrature can be used to evaluate the surface integrals and this results in about 60\% savings in coding effort. The new formulation also affords a numerically robust way to calculate the potential close to the boundaries. Detailed implementations of this approach are illustrated with problems involving osculating boundaries, 2D domains with corners and a wave drag problem in a 3D semi-infinite domain. The explicit formulation of problems with axial symmetry is also given.
\end{abstract}


\end{frontmatter}

% \linenumbers


% ---- Section 1  Introduction ----
\section{Introduction}
\label{}
Solution of the Laplace equation for the potential problem underpins many applications in electrostatics and heat conduction. It is also central to modelling moving or deformable boundaries in the high Reynolds number regime in fluid mechanics. There, viscous and boundary layer effects are not dominant whereby a description based on potential flow can therefore be used as a first approximation or as a base case for further refinement. Numerous examples can be found in hydrofoil dynamics~\citep{Faltinsen2008, Xu2013a}, the description of waves~\citep{Xue2001, Liu2001}, cavitation or supercavitation phenomena~\citep{Blake1986, Blake1987}, in civil, marine and ocean engineering and oscillating bubble dynamics in sonophysics and sonochemistry~\citep{Leighton1994}.

Many of the above applications in multiphase fluid mechanics require the accurate tracking of moving interfaces that can be cumbersome and expensive to implement using grid based methods, especially in 3D. Consequently the use of the boundary element method (BEM) is advantageous because computational effort can focus on modelling all the important interfaces with the additional benefit of reducing the dimension of the problem by one, thus obviating the need to compute solutions in the whole flow domain. To track moving and deforming interfaces with precision, the Laplace equation is solved at each time step and the boundaries evolve according to the unsteady Bernoulli equation. This approach can readily be adapted to handle large or infinite domains or boundaries. The BEM is especially appealing for infinite fluid domains since the behaviour at infinity can be accounted for analytically. Although the BEM generates dense matrix equations, the Fast Multipole Boundary Element Method~\citep{Liu2009} can be used to reduce both CPU time and memory requirement from $O(N^2)$ to $O(N \log N)$. Thus in spite of being a simplification of the full Navier-Stokes description, the theory of potential flow together with the Bernoulli equation to describe unsteady problems occupies an important role in multiphase fluid dynamics.

The inherent use of the fundamental solution in the formulation of the BEM means that the integral equation contains singular kernels. This characteristic feature has been described as ``a mathematical monster that leaps out of every page'' due to ``very unfamiliar and complex mathematics''~\citep{Becker92}. Since the physical problem itself is perfectly well behaved on the boundaries, such singularities are numerical inconveniences generated by the mathematical formulation~\citep{LiuRudolphi1999}. Traditionally, the singular behaviour is dealt with by a local change of variables in the evaluation of the surface integrals~\citep{Telles1987} that comes at the expense of additional coding effort. Previous attempts to remove such singularities analytically required the introduction of additional unknowns such as derivatives tangential to the surface that have to be found by developing and solving extra integral equations~\citep{LiuRudolphi1999, Chen2005}. Another method to remove the singularities requires finding additional parameters that have to be determined on a fictitious `nearby' boundary~\citep{Cao1991, Zhang1999, Chen2009}.

Moreover, when two different boundaries or two parts of the same boundary become close to each other, the traditional implementation of the BIM does not prevent the deleterious influence of singularities that originate from nodes of one boundary on the other. In particular, for moving boundary potential problems, it is highly desirable to eliminate all singular terms that arise in traditional formulations of the BEM as this avoids the need to track the spatial separation of different parts of the boundary and to determine when remedial action may be required.

Here we show that the well-known mathematical singularities that arise in the BEM for solving the Laplace equation can be removed analytically \emph{without} generating additional unknowns or equations to be solved. The desingularised formulation is given for both general 3D and axisymmetric cases. The approach can also be applied to evaluate the potential at points near boundaries in a numerically robust way. The implementation and resulting improvement in accuracy are illustrated with a number of examples: a problem with osculating boundaries, a 3D problem with a semi-infinite domain that arises in the study of wave drag near a deformable surface (with movies in the electronic supplement) and problems involving domains with corners in 2D.

% ---- Section 2  NS BIM for general psi ----
\section{Non-singular formulation of the boundary integral method}
\begin{figure}[t]
\centering
\includegraphics[width=0.6\textwidth] {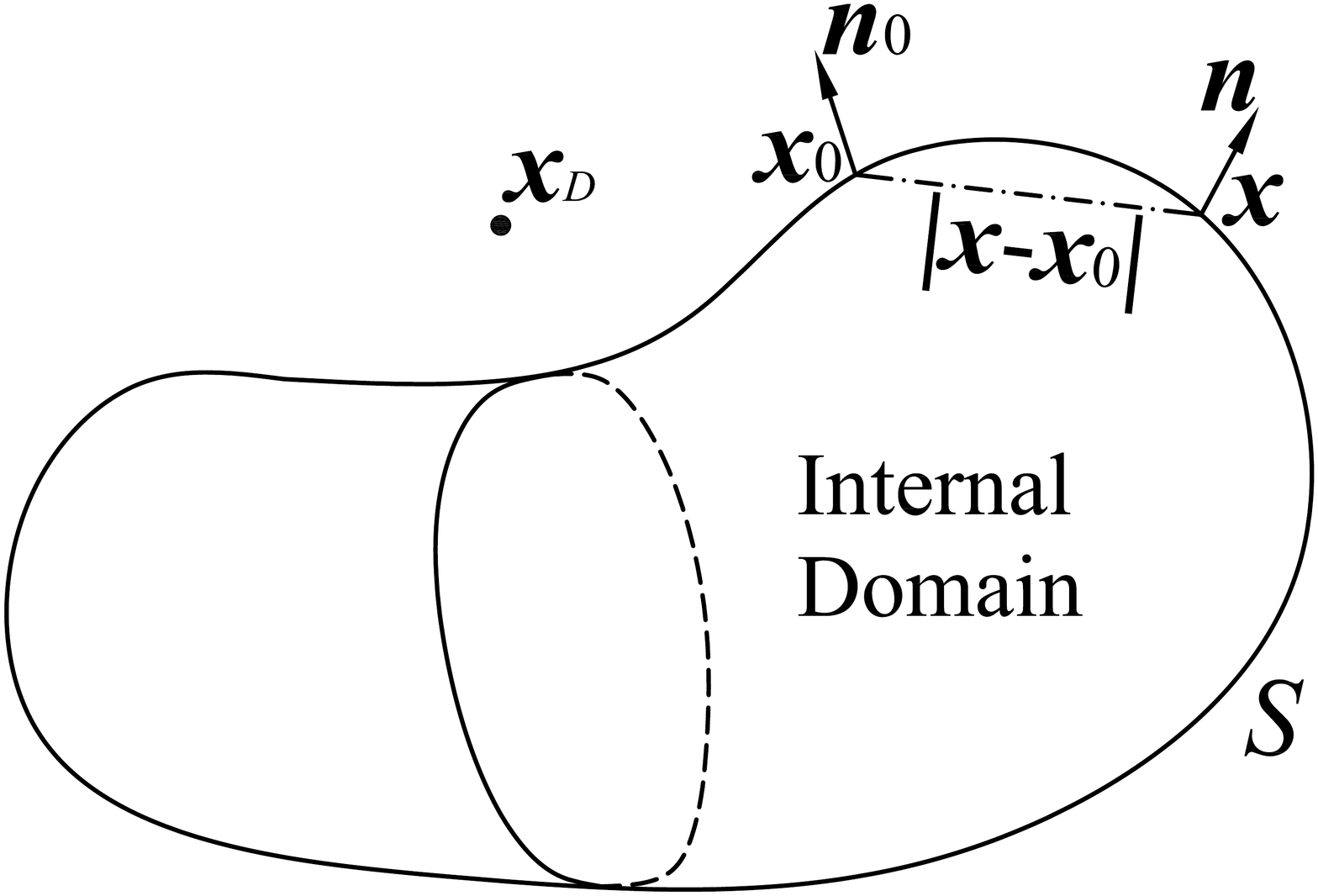}
\caption{The 3D internal domain defined by the closed surface $S$ with the observation point $\boldsymbol{x_0}$ with outward normal $\boldsymbol{n}_{0}$, the integration point $\boldsymbol{x}$ with outward normal $\boldsymbol{n}$ and a sample location of $\boldsymbol{x}_{D}$.}\label{fig:illus_internal}
\end{figure}

To develop a non-singular formulation of the boundary integral method, consider the internal  problem in a 3D domain that is bounded by the closed surface $S$, as shown in Fig.~\ref{fig:illus_internal}. The potential, $\phi$, is governed by the Laplace equation
\begin{align}\label{eq:lplc}
\bnabla^2\phi = 0.
\end{align} %that follows from the continuity condition: $\bnabla\bcdot\boldsymbol{u}=0$.
By using the 3D free space Green's function $G(\boldsymbol{x},\boldsymbol{x}_0) = 1/|\boldsymbol{x}-\boldsymbol{x}_0|$ and with the help of Green's second identity, the solution of Eq.~(\ref{eq:lplc}) can be found by solving the conventional boundary integral equation
\begin{align}\label{eq:cnvtnbim}
c_{0}\phi(\boldsymbol{x}_0) + \int_{S} \phi(\boldsymbol{x}) \frac{\partial{G}}{\partial{n}} \; \text{d}S(\boldsymbol{x}) = \int_{S} \frac{\partial{\phi}}{\partial{n}} G(\boldsymbol{x},\boldsymbol{x}_0) \; \text{d}S(\boldsymbol{x}),
\end{align}
where $c_0$ is the solid angle at the observation point $\boldsymbol{x}_0$ on $S$ and the surface element $\text{d}S(\boldsymbol{x})$ is at $\boldsymbol{x}$. The normal derivatives are defined by
${\partial{\phi}}/{\partial{n}} \equiv \bnabla\phi(\boldsymbol{x})\bcdot\boldsymbol{n(\boldsymbol{x})}$
and
${\p{G}}/{\p{n}} \equiv \bnabla G(\boldsymbol{x},\boldsymbol{x}_0)\bcdot\boldsymbol{n(\boldsymbol{x})}$,
where $\boldsymbol{n} \equiv \boldsymbol{n(\boldsymbol{x})}$ is the unit normal vector pointing out of the internal domain at $\boldsymbol{x}$~\citep{Becker92}. If either $\phi$ or ${\p{\phi}}/{\p{n}}$ or a mixed boundary condition is given on the whole boundary $S$, the corresponding ${\p{\phi}}/{\p{n}}$ or/and $\phi$ can be obtained from Eq.~(\ref{eq:cnvtnbim}). Although numerical integration over the singularities in $G(\boldsymbol{x},\boldsymbol{x}_0)$ and ${\partial{G}}/{\partial{n}}$ can be effected by established methods~\citep{Telles1987}, our aim is to remove such mathematical singularities analytically at the outset.

% -- FIGURE 1 --

Corresponding to a given point $\boldsymbol{x}_0$ on the boundary, we construct a function $\psi(\boldsymbol{x})$ that also satisfies the Laplace equation and hence Eq. (\ref{eq:cnvtnbim}), with the properties that as $\boldsymbol{x} \rightarrow \boldsymbol{x}_0$, $\psi(\boldsymbol{x}) \rightarrow \phi(\boldsymbol{x}_0) \equiv \phi_0$ and ${\partial{\psi}}/{\partial{n}} \rightarrow ({\partial{\phi}}/{\partial{n}})_0 \equiv \bnabla \phi(\boldsymbol{x}_0) \bcdot \boldsymbol{n}(\boldsymbol{x}_0)$, where $\boldsymbol{n}_{0} \equiv \boldsymbol{n}(\boldsymbol{x}_0)$ is the outward unit normal vector at $\boldsymbol{x}_0$. We choose $\psi(\boldsymbol{x})$ to be of the form
\begin{align}\label{eq:psi_f(x)}
\psi(\boldsymbol{x}) \equiv \phi(\boldsymbol{x}_0) + \left(\frac{\partial{\phi}}{\partial{n}}\right)_{0} f(\boldsymbol{x})
\end{align}
so the function $f(\boldsymbol{x})$ must satisfy
\begin{align}\label{eq:f(x)}
\bnabla^2 f(\boldsymbol{x}) = 0, \qquad f(\boldsymbol{x}_0) = 0, \qquad \bnabla f(\boldsymbol{x}_0) \bcdot \boldsymbol{n}_0 = 1.
\end{align}
Taking the difference between the conventional boundary integral equations for $\phi(\boldsymbol{x})$ and for $\psi(\boldsymbol{x})$ we obtain an integral equation relating $\phi(\boldsymbol{x})$ and $\partial{\phi}/\partial{n}$ on $S$ that replaces the conventional boundary integral equation in Eq.~(\ref{eq:cnvtnbim}):
\begin{align}\label{eq:nsbim_f}
\int_{S} \left[\phi(\boldsymbol{x}) - \phi(\boldsymbol{x}_0) - \left(\frac{\partial{\phi}}{\partial{n}}\right)_{0} f(\boldsymbol{x})\right]
 \frac{\partial{G(\boldsymbol{x},\boldsymbol{x}_0)}}{\partial{n}} \; \text{d}S(\boldsymbol{x}) \qquad \qquad \qquad \nonumber \\
= \int_{S} \left[ \frac{\partial{\phi}}{\partial{n}} - \left(\frac{\partial{\phi}}{\partial{n}}\right)_{0} \bnabla f(\boldsymbol{x}) \bcdot \boldsymbol{n} \right] G(\boldsymbol{x},\boldsymbol{x}_0) \; \text{d}S(\boldsymbol{x}).
\end{align}
The key point is that both integrands in Eq.~(\ref{eq:nsbim_f}) are now non-singular and thus any convenient quadrature method can be used to evaluate the integrals. Note also that the solid angle $c_{0}$ no longer appears. Implicit in the derivation of Eq.~(\ref{eq:nsbim_f}) is that the outward unit normal $\boldsymbol{n}(\boldsymbol{x}_0)$ is uniquely defined at $\boldsymbol{x}_0$. The implementation of Eq.~(\ref{eq:nsbim_f}) at nodes where the normal is not defined, e.g. at a corner, is considered in Section 5. A detailed proof of this de-singularisation method using the linear function: $f(\boldsymbol{x}) = \boldsymbol{n}_0 \bcdot (\boldsymbol{x} - \boldsymbol{x}_0)$ has been given elsewhere~\cite{Klaseboer2012} and establishes the theoretical basis of the numerical scheme constructed earlier to regularise the system of linear equations that arise from the standard implementation of the BEM~\citep{Klaseboer2009}. This approach has also been extended to de-singularise boundary integral equations that arise in Stokes flow, in solving the Helmholtz equation and the equations of linear elasticity using a linear function for $f(\boldsymbol{x})$~\citep{Klaseboer2012}.

However, for potential problems involving domains of infinite extent (external problems), a different choice of  $f(\boldsymbol{x})$ is required since the linear function is unbounded at infinity. One possible choice for $f(\boldsymbol{x})$ is
\begin{eqnarray}\label{eq:psi_inverse}
\psi(\boldsymbol{x})
&\equiv& \phi(\boldsymbol{x}_0) + \left(\frac{\partial{\phi}}{\partial{n}}\right)_{0} f(\boldsymbol{x}) \nonumber \\
&=& \phi(\boldsymbol{x}_0) + \left(\frac{\partial{\phi}}{\partial{n}}\right)_{0} \frac{|\boldsymbol{x}_0 - \boldsymbol{x}_D|^2}{\boldsymbol{n}_0 \bcdot (\boldsymbol{x}_0 - \boldsymbol{x}_D)} \left( 1-\frac{|\boldsymbol{x}_0 - \boldsymbol{x}_D|}{|\boldsymbol{x} - \boldsymbol{x}_D|}\right).
\end{eqnarray}
The constant vector $\boldsymbol{x}_D$ is the position of any convenient point that is located \emph{outside} the solution domain and satisfies $\boldsymbol{n}_0 \bcdot (\boldsymbol{x}_0 - \boldsymbol{x}_D) \neq 0$. In the next section, we will use this form of $\psi(\boldsymbol{x})$ to formulate a non-singular BEM for axisymmetric problems and in Section 6, we demonstrate how Eq.~(\ref{eq:psi_inverse}) can be used to formulate a non-singular boundary integral problem in a semi-infinite domain to solve a 3D time-dependent potential flow problem. Before proceeding, we note that we are not restricted to using a linear function or the form given by Eq.~(\ref{eq:psi_inverse}) for $f(\boldsymbol{x})$ to construct non-singular versions of the boundary integral equation. In fact, any form of $f(\boldsymbol{x})$ that satisfies the conditions given by Eq.~(\ref{eq:f(x)}) can be considered for use in Eq.~(\ref{eq:nsbim_f}) to remove the singular behaviour due to the presence of $G$ and ${\partial{G}}/{\partial{n}}$ at $\boldsymbol{x} =\boldsymbol{x}_0$.

% ---- Section 3  Axisymmetric NS BIM ----
\section{Non-singular axisymmetric boundary integral equation}

For problems that possess axial symmetry whereby in the cylindrical variables: $r, \theta$ and $z$, we have $\phi(\boldsymbol{x}) = \phi(r,z)$ and $\partial{\phi(\boldsymbol{x})}/\partial{n} = \partial{\phi(r, z)}/\partial{n}$, the integration over the azimuthal angle, $\theta$, can be evaluated analytically~\citep{Becker92}. Since a point $(r,z)$ on the axisymmetric boundary surface is specified by some given equation $S(r, z) = 0$, the surface integrals for the axisymmetric case can be reduced to 1D integrals.

The $\psi(\boldsymbol{x})$ given by Eq.~(\ref{eq:psi_inverse}) will be an axisymmetric function if we choose (in Cartesian coordinates) $\boldsymbol{x}_D = (0,0,z_D)$ to lie on the $z$-axis of symmetry with $z_D$ located outside the solution domain. We measure the azimuthal angle, $\theta$, relative to $\boldsymbol{x}_0 = (r_0,0,z_0)$ so that $\boldsymbol{x} = (r \cos \theta,r \sin \theta,z)$ and the surface normals are given by $\boldsymbol{n}_0 = (n_{r0},0,n_{z0})$ and $\boldsymbol{n} = (n_{r} \cos \theta,n_{r} \sin \theta,n_{z})$. With these definitions, the axisymmetric function $\psi(\boldsymbol{x})$ in Eq. (\ref{eq:psi_inverse}) and its normal derivative are given explicitly by
\begin{eqnarray}
  \psi(r,z) &=& \phi(r_0,z_0) + \left(\frac{\partial{\phi}}{\partial{n}}\right)_{0} \left ( \frac{\rho -\rho_0}{\rho}  \right ) \frac{\rho_0^2}{s_0} \label{eq:psi_axisym} \\
 \nonumber \\
  \frac{\partial{\psi(r,z)}}{\partial{n}} &=& \left(\frac{\partial{\phi}}{\partial{n}}\right)_{0} \left ( \frac{\rho_0}{\rho} \right )^3 \left ( \frac{s}{s_0} \right ) \label{eq:d_psi_dn_axisym}
\end{eqnarray}
where $\rho \equiv \sqrt{r^2 +(z - z_D)^2}$, $\rho_0 \equiv \sqrt{r_{0}^2 +(z_0 - z_D)^2}$, $s \equiv r n_r + (z - z_D) n_z$ and $s_0 \equiv r_0 n_{r0} + (z_0 - z_D) n_{z0}$.

The axisymmetric version of the non-singular boundary integral equation obtained after performing the $\theta$-integration in Eq.~(\ref{eq:nsbim_f}) using Eq.~(\ref{eq:psi_inverse}) is (see also~\citep{Becker92})
\begin{eqnarray}\label{eq:nsbim_axisym}
\int \frac{\chi(r,z) \, r \, E(m)}{(1-m) \bar{R}^3} \left[ (r - r_0)n_r + (z - z_0)n_z - 2 \, r_0 \, n_r \, (1-m)/m \right] \text{d}{\Gamma} \nonumber \\
=\int \frac{2 \, r \, r_0 \, n_r \, \chi(r,z) \, K(m)}{m \bar{R}^3} \; \text{d}{\Gamma}  +\int  \frac{\partial{\chi(r,z)}}{\partial{n}} \, \frac{r \, K(m)}{\bar{R}} \; \text{d}{\Gamma}
\end{eqnarray}
where $\chi(r,z) \equiv \phi(r,z) - \psi(r,z)$. The arc length element $\text{d}{\Gamma}$ may be expressed as $\sqrt{1+(\text{d}r/\text{d}z)^2}\text{ d}z$ using the equation: $S(r,z) = 0$ that defines the axisymmetric surface. The quantities $K(m)$ and $E(m)$ are the complete elliptic integrals of the first and second kind of the parameter, $m$ (see~\citep{Becker92, Abramowitz1972}), with
\begin{eqnarray}
  m \equiv \frac{4 r r_o}{(r + r_0)^2 + (z - z_0)^2}\equiv \frac{4 r r_o}{\bar{R}^2},  \quad 0 \le m \le 1.
\end{eqnarray}
The axisymmetric boundary integral equation given by Eq.~(\ref{eq:nsbim_axisym}) is non-singular because the $\log(1-m)$ divergence in $K(m)$ as $m \rightarrow 1$ in the limit $\boldsymbol{x} \rightarrow \boldsymbol{x}_0$ is now suppressed by the terms containing the difference between $\phi$ and $\psi$ and between ${\partial{\phi}}/{\partial{n}}$ and ${\partial{\psi}}/{\partial{n}}$ that vanish as $(1-m)$ as $\boldsymbol{x} \rightarrow \boldsymbol{x}_0$. Thus the line integral in Eq.~(\ref{eq:nsbim_axisym}) can be evaluated using any quadrature method. This is a consequence of using an axisymmetric form of $\psi(r,z)$ given by Eq.~(\ref{eq:psi_inverse}) that simplified to Eqs.~(\ref{eq:psi_axisym}) and~(\ref{eq:d_psi_dn_axisym}).

An additional bonus of the non-singular axisymmetric formulation in Eq. (\ref{eq:nsbim_axisym}) is that no special effort is required to handle node points that are located on the axis of symmetry that has been a technical inconvenience of the conventional axisymmetric form of the boundary integral equation~\citep{Becker92}.

% ---- Section 4  Potential evaluation ----
\section{Robust method to calculate the potential near a boundary}

The method of de-singularising the boundary integral equation by subtracting out the singular behaviour can also be used to give a numerically robust method to calculate the value of the potential at points near boundaries. To find the potential, $\phi(\boldsymbol{x}_p)$, at an observation point $\boldsymbol{x}_p$ that is located inside the domain but may be close to the boundary, we begin with the conventional boundary equation, Eq.~(\ref{eq:cnvtnbim}), for the function $[\phi(\boldsymbol{x}_p) - \psi(\boldsymbol{x}_p)]$, with $\psi(\boldsymbol{x})$ given by Eq.~(\ref{eq:psi_f(x)}) and $c_{0}=4\pi$ for points inside the domain,
\begin{align} \label{eq:conventional_phi_psi}
  4 \pi \left\{ \phi(\boldsymbol{x}_p) - \left[ \phi(\boldsymbol{x}_0) +   \left( \frac{\partial{\phi}} {\partial{n}} \right)_{0} f(\boldsymbol{x}_p) \right] \right\} \qquad \qquad \qquad \qquad \qquad  &&
   \nonumber \\
   \qquad \qquad  + \int \left\{ \phi(\boldsymbol{x}) - \left[ \phi(\boldsymbol{x}_0) + \left( \frac{\partial{\phi}}{\partial{n}} \right)_0 f(\boldsymbol{x}) \right]  \right\} \frac {\partial{G(\boldsymbol{x}, \boldsymbol{x}_p)}} {\partial{n}} \; \text{d}S(\boldsymbol{x}) &&
  \nonumber \\
  = \int \left\{ \frac{\partial{\phi}}{\partial{n}} - \left( \frac{\partial{\phi}}{\partial{n}} \right)_0 \bnabla f(\boldsymbol{x}) \bcdot \boldsymbol{n}(\boldsymbol{x}) \right\} G(\boldsymbol{x}, \boldsymbol{x}_p) \; \text{d}S(\boldsymbol{x}).
\end{align}
The point $\boldsymbol{x}_0$ is located on the boundary and its relation to $\boldsymbol{x}_p$ is specified below. The nearly singular behaviour of the integrands when $\boldsymbol{x}_p$ is close to the boundary can be eliminated by subtracting Eq.~(\ref{eq:nsbim_f}) from Eq.~(\ref{eq:conventional_phi_psi}) to give
\begin{align} \label{eq:nonsig_phi_xp_f}
  \phi(\boldsymbol{x}_p) =  \phi(\boldsymbol{x}_0) +   \left( \frac{\partial{\phi}} {\partial{n}} \right)_{0} f(\boldsymbol{x}_p)  \qquad \qquad \qquad \qquad \qquad \qquad \qquad \qquad &&
   \nonumber \\
  - \frac {1}{4 \pi} \int \left\{ \phi(\boldsymbol{x}) - \left[ \phi(\boldsymbol{x}_0) + \left( \frac{\partial{\phi}}{\partial{n}} \right)_0 f(\boldsymbol{x}) \right] \right\} \left\{ \frac {\partial{G(\boldsymbol{x}, \boldsymbol{x}_p)}} {\partial{n}} -   \frac {\partial{G(\boldsymbol{x}, \boldsymbol{x}_0)}} {\partial{n}} \right\} \; \text{d}S(\boldsymbol{x}) &&
  \nonumber \\
  + \frac {1}{4 \pi}  \int \left\{ \frac{\partial{\phi}}{\partial{n}} - \left( \frac{\partial{\phi}}{\partial{n}} \right)_0 \bnabla f(\boldsymbol{x}) \bcdot \boldsymbol{n}(\boldsymbol{x}) \right\} \left\{ G(\boldsymbol{x}, \boldsymbol{x}_p) - G(\boldsymbol{x}, \boldsymbol{x}_0) \right\} \; \text{d}S(\boldsymbol{x}).
\end{align}
When the observation point $\boldsymbol{x}_p$ is near the boundary, a suitable choice for the boundary point, $\boldsymbol{x}_0$, is to project $\boldsymbol{x}_p$ onto the boundary using the relation $\boldsymbol{x}_p = \boldsymbol{x}_0 - \varepsilon \boldsymbol{n}_0$, with $\varepsilon > 0$.

Another important benefit arising from the numerical robustness of the present de-singualarised formulation of the BIM occurs when different parts of the domain boundary become close together as often happens in multi-scale moving boundary problems. A simple example of this is when the boundary comprises of two nearly touching spheres and the point $\boldsymbol{x}_0$ is located near where the spherical boundaries are close together. The present formulation eliminates the near singular behaviour of the kernel at nodes that are on a different part of boundary that is close spatially to the observation point $\boldsymbol{x}_0$. A numerical demonstration of this is given in Section 6.

% ---- Section 5 - Corner nodes ----
\section{Corner nodes in 2D}
The above formulation of the non-singular form of the BEM for the potential problem assumes that the unit normal $n_0 \equiv n(\boldsymbol{x}_0)$ is uniquely defined at $\boldsymbol{x}_0$. When $\boldsymbol{x}_0$ is a corner node on a 2D boundary, we choose $\psi(\boldsymbol{x})$ to have the form
\begin{align}\label{eq:crpsi_f(x)}
\psi(\boldsymbol{x}) = \phi(\boldsymbol{x}_0) + \left(\frac{\partial{\phi}}{\partial{n}}\right)^{L}_{0} f^{L}(\boldsymbol{x}) + \left(\frac{\partial{\phi}}{\partial{n}}\right)^{R}_{0} f^{R}(\boldsymbol{x})
\end{align}
where $\phi_{0} \equiv \phi(\boldsymbol{x}_0)$ is the potential at $\boldsymbol{x}_0$, and the constants $(\p{\phi}/\p{n})^{L}_{0}$ and $(\p{\phi}/\p{n})^{R}_{0}$ are the values of normal derivatives as one approaches $\boldsymbol{x}_0$ from the `left' ($L$) and from the `right' ($R$). The functions $f^{L}(\boldsymbol{x})$ and $f^{R}(\boldsymbol{x})$ satisfy the following conditions
\begin{align}\label{eq:crf(x)L}
\bnabla^2 f^{L}(\boldsymbol{x}) = 0, \quad f^{L}(\boldsymbol{x}_0) = 0, \quad \bnabla f^{L}(\boldsymbol{x}_0) \bcdot \boldsymbol{n}^{L}_0 = 1, \quad \bnabla f^{L}(\boldsymbol{x}_0) \bcdot \boldsymbol{n}^{R}_0 = 0.
\end{align}
\begin{align}\label{eq:crf(x)R}
\bnabla^2 f^{R}(\boldsymbol{x}) = 0, \quad f^{R}(\boldsymbol{x}_0) = 0, \quad \bnabla f^{R}(\boldsymbol{x}_0) \bcdot \boldsymbol{n}^{R}_0 = 1, \quad \bnabla f^{R}(\boldsymbol{x}_0) \bcdot \boldsymbol{n}^{L}_0 = 0.
\end{align}

The value of $\nabla \phi(\boldsymbol{x}_0)$ at the corner $\boldsymbol{x}_0$ can be expressed in two equivalent forms
\begin{align}\label{eq:compcon}
\left\{
\begin{aligned}
&\nabla \phi(\boldsymbol{x}_0) = \left(\frac{\partial{\phi}}{\partial{n}}\right)^{L}_{0}\boldsymbol{n}^{L}_{0} + \left(\frac{\partial{\phi}}{\partial{t}}\right)^{L}_{0}\boldsymbol{t}^{L}_{0}, \\
&\nabla \phi(\boldsymbol{x}_0) = \left(\frac{\partial{\phi}}{\partial{n}}\right)^{R}_{0}\boldsymbol{n}^{R}_{0} + \left(\frac{\partial{\phi}}{\partial{t}}\right)^{R}_{0}\boldsymbol{t}^{R}_{0};
\end{aligned}
\right.
\end{align}
in which, $(\partial{\phi}/\partial{t})^{L}_{0}$ and $(\partial{\phi}/\partial{t})^{R}_{0}$ are the two tangential derivatives of the potential on side $L$ and $R$ at $\boldsymbol{x}_{0}$, respectively, and $\boldsymbol{t}^{L}_{0}$ and $\boldsymbol{t}^{R}_{0}$ are the unit tangential vectors along side $L$ and $R$ at $\boldsymbol{x}_{0}$, respectively. This compatibility condition in Eq.~(\ref{eq:compcon}) provides an additional relation between $(\partial{\phi}/\partial{n})^{L}_{0}$ and $(\partial{\phi}/\partial{n})^{R}_{0}$ ~\citep{Grilli1990} in the formulation of the BIM whereas the tangential derivatives can be constructed from the values of $\phi$ at neighbouring nodes.

A choice for $f^{L}(\boldsymbol{x})$ and $f^{R}(\boldsymbol{x})$ for a 2D corner problem can be constructed in terms of the surface normals $\boldsymbol{n}^{L}_0$ and $\boldsymbol{n}^{R}_0$ on either side of the corner at $\boldsymbol{x}_0$
\begin{align}\label{eq:crnrchc1L}
f^{L}(\boldsymbol{x}) = \frac{-(\boldsymbol{n}^{L}_0 \bcdot \boldsymbol{n}^{R}_0)\boldsymbol{n}^{R}_0 + \boldsymbol{n}^{L}_0} {1-(\boldsymbol{n}^{L}_0 \bcdot \boldsymbol{n}^{R}_0)^2} \bcdot (\boldsymbol{x}-\boldsymbol{x}_0),
\end{align}
and
\begin{align}\label{eq:crnrchc1R}
f^{R}(\boldsymbol{x}) = \frac{-(\boldsymbol{n}^{R}_0 \bcdot \boldsymbol{n}^{L}_0)\boldsymbol{n}^{L}_0 + \boldsymbol{n}^{R}_0} {1-(\boldsymbol{n}^{R}_0 \bcdot \boldsymbol{n}^{L}_0)^2} \bcdot (\boldsymbol{x}-\boldsymbol{x}_0).
\end{align}

% ---- Section 6 - Numerical examples ----
\section{Numerical demonstrations and examples}

% -- FIGURE 2 --
\begin{figure}[t]
\centering
\subfloat[]{ \includegraphics[width=0.45\textwidth] {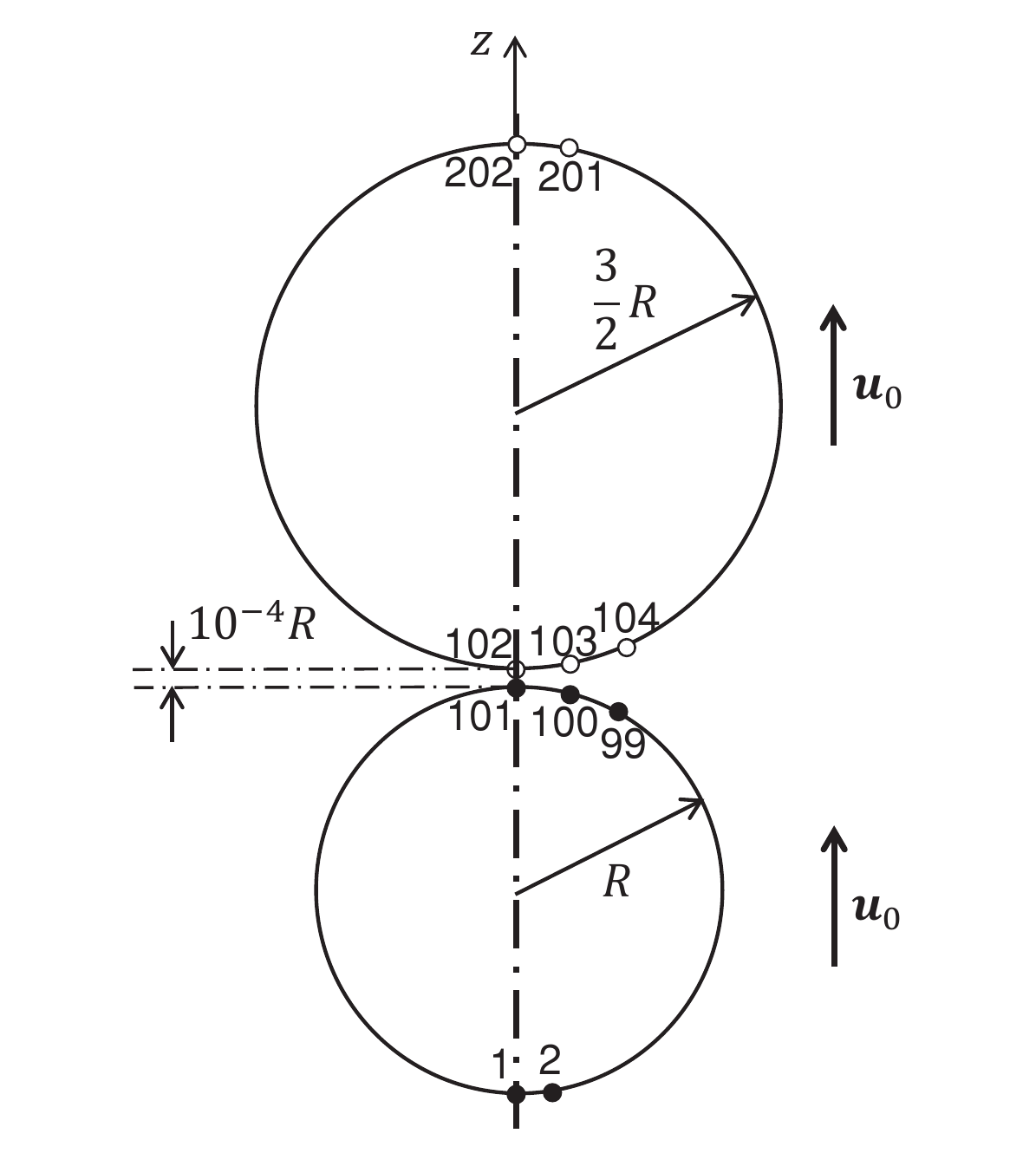} }
\subfloat[]{ \includegraphics[width=0.52\textwidth] {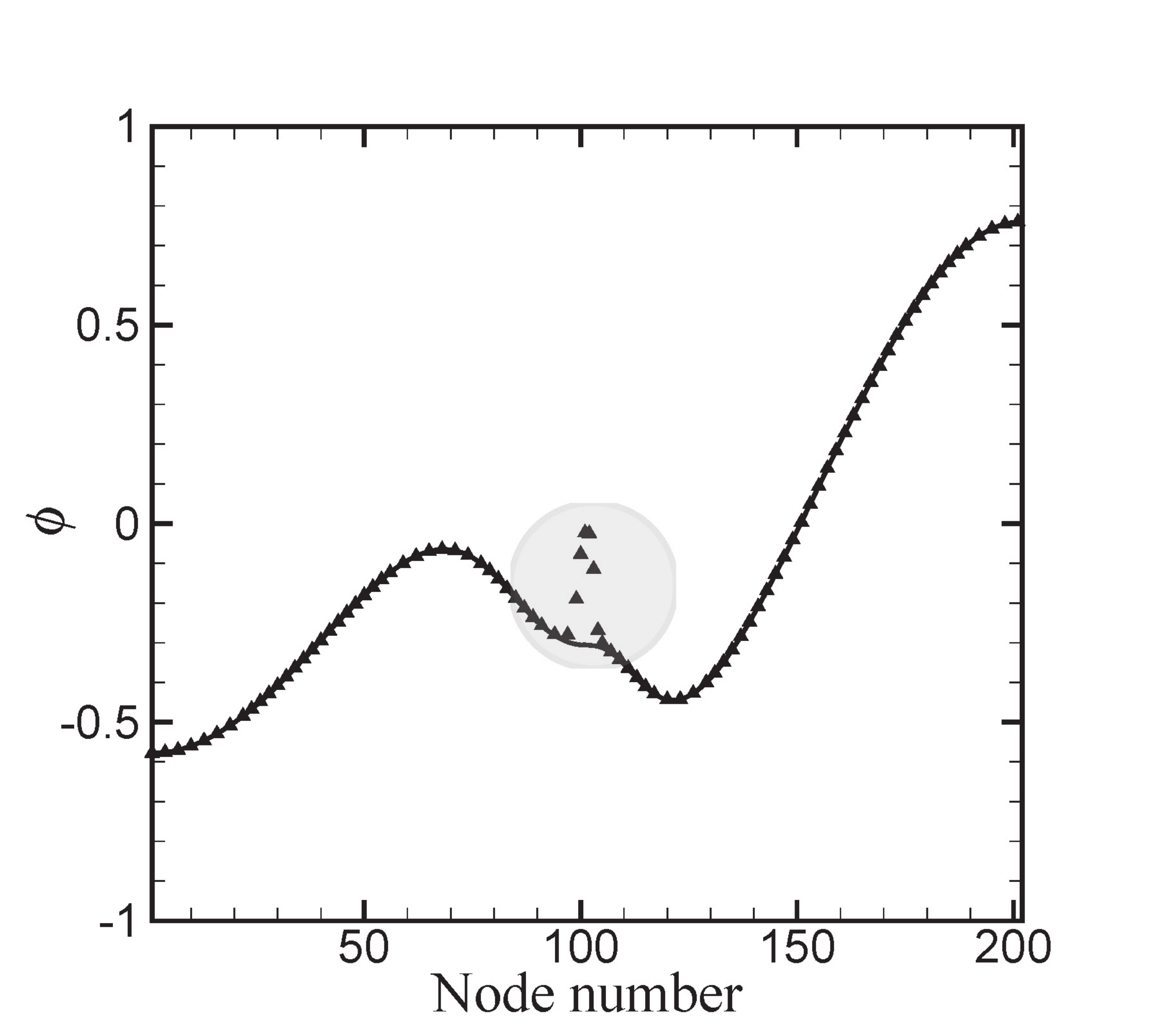} }
\caption{Variation of the velocity potential along the meridian for two nearly osculating spheres translating along their line of centre at equal constant speed in potential flow. (a) Nodes 1--101 are on the sphere of radius $R$ and nodes 102--202 are on the sphere of radius $3R/2$. The gap between the spheres is $10^{-4} R$. (b) Erroneous results of the potential obtained from the standard axisymmetric boundary integral method (symbols) when compared to our non-singular method given by Eq.~(\ref{eq:nsbim_axisym}) (line) are circled.}\label{fig:miloh3}
\end{figure}

\subsection{Two nearly touching spheres}
To demonstrate the utility and robustness of our non-singular boundary integral method for problems with near osculating boundaries, we consider the potential problem associated with two nearly touching spheres translating along their line of centre at identical constant speed. The spheres have radii $R$ and $3R/2$ and are at a separation of $10^{-4} R$ at the point of closest approach (Fig.~\ref{fig:miloh3}a). We compare the velocity potential obtained using the standard axisymmetric boundary integral method~\citep{Wang1996}, with our non-singular version given by Eq.~(\ref{eq:nsbim_axisym}) with $\boldsymbol{x}_{D}$ set to the centre of each sphere. In Fig.~\ref{fig:miloh3}b, we see that the error in the standard axisymmetric boundary integral method is very large in the region where the spheres are close together. The standard 3D boundary integral method gives errors similar to the standard axisymmetric method whereas the non-singular 3D version given by Eq.~(\ref{eq:nsbim_f}) has the same accuracy as the axisymmetric version. This large error arises in the standard version of the boundary integral method because of the influence of the observation point on one sphere from the singular kernel centred at nearby nodes located on the other sphere.  As we have shown in Section 4, with our non-singular formulation, such effects do not arise.

\subsection{Corner problem in 2D}

% -- FIGURE 3 --
\begin{figure}[t]
\centering
\includegraphics[width=0.6\textwidth] {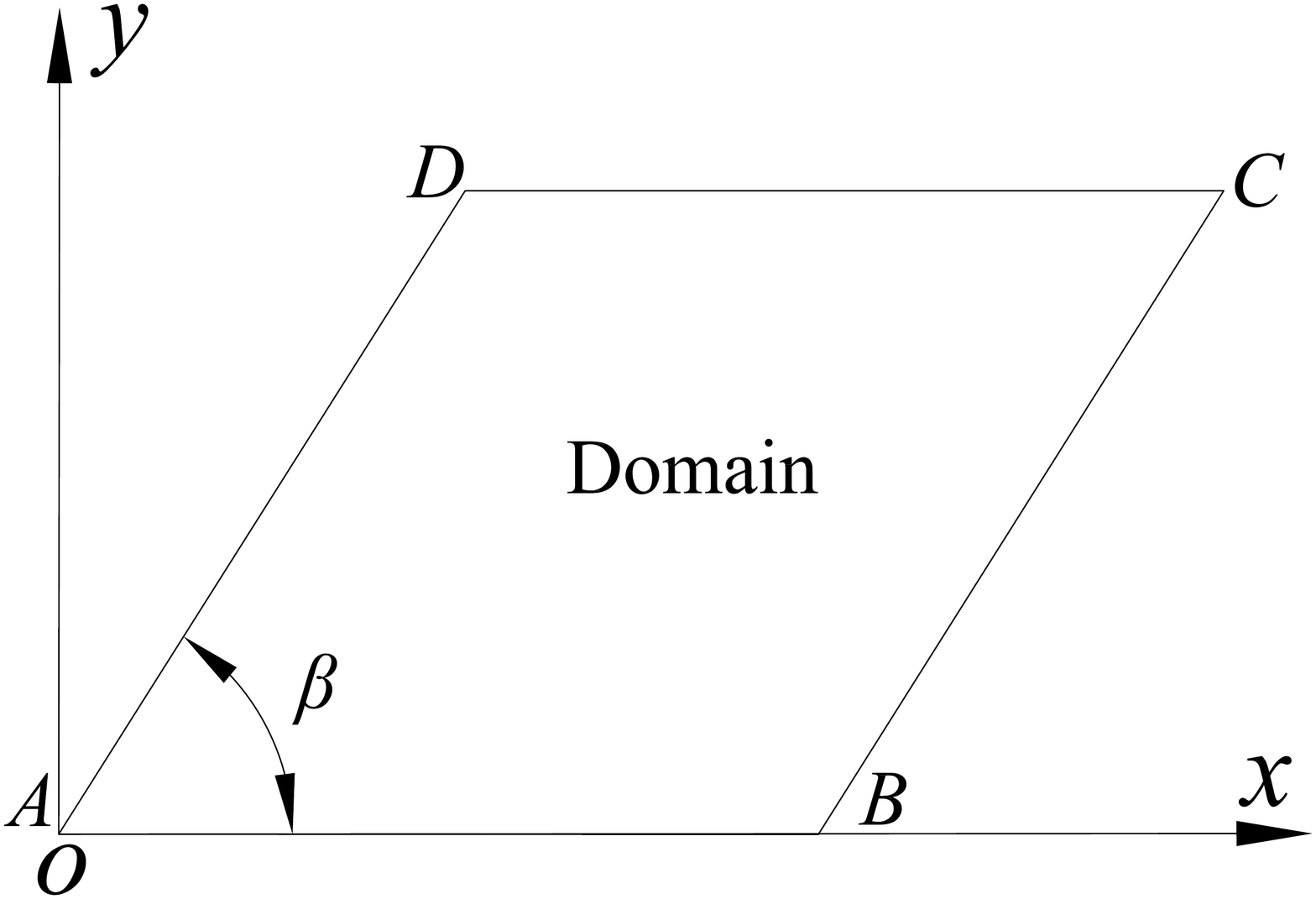}
\caption{Rectangular domain for an interior 2D potential problems.}\label{fig:cornersketch}
\end{figure}

To illustrate the implementation of our non-singular formulation of the BIM, we consider three different interior potential problems in rectangular domains of different shapes, as shown in Fig.~\ref{fig:cornersketch}. The length of all four edges are set to be 1 unit. Uniform linear elements are employed on the boundary of the parallelogram specified by the angle, $\beta$. At the corner nodes, the function $\psi(\boldsymbol{x})$ given by Eq.~(\ref{eq:crpsi_f(x)}) used in the non-singular boundary integral equation~(\ref{eq:nsbim_f}) and the double node technique~\cite{Grilli1990} has been applied. The tangential derivatives of the potential in the compatibility condition in Eq. (\ref{eq:compcon}) are obtained by a fourth order finite difference scheme using values of $\phi$ along the boundaries. For the remaining nodes on the edges, the form of $\psi(\boldsymbol{x})$ in Eq.~(\ref{eq:psi_f(x)}) with $f(\boldsymbol{x}) = \boldsymbol{n}_0 \bcdot (\boldsymbol{x} - \boldsymbol{x}_0)$ is used in Eq.~(\ref{eq:nsbim_f}).

We tested three cases corresponding to analytical solutions:
\begin{align}
\phi = 1-x, \quad \text{Case I};
\end{align}
\begin{align}
\phi = 1-xy, \quad \text{Case II};
\end{align}
\begin{align}
\phi = \frac{\sinh(\pi y)}{\sinh(\pi)}\sin(\pi x) + 1, \quad \text{Case III}.
\end{align}

% -- FIGURE 4 --
\begin{figure}[t]
\centering
\subfloat[$\beta = 45^{\text{o}}$]{ \includegraphics[width=0.5\textwidth] {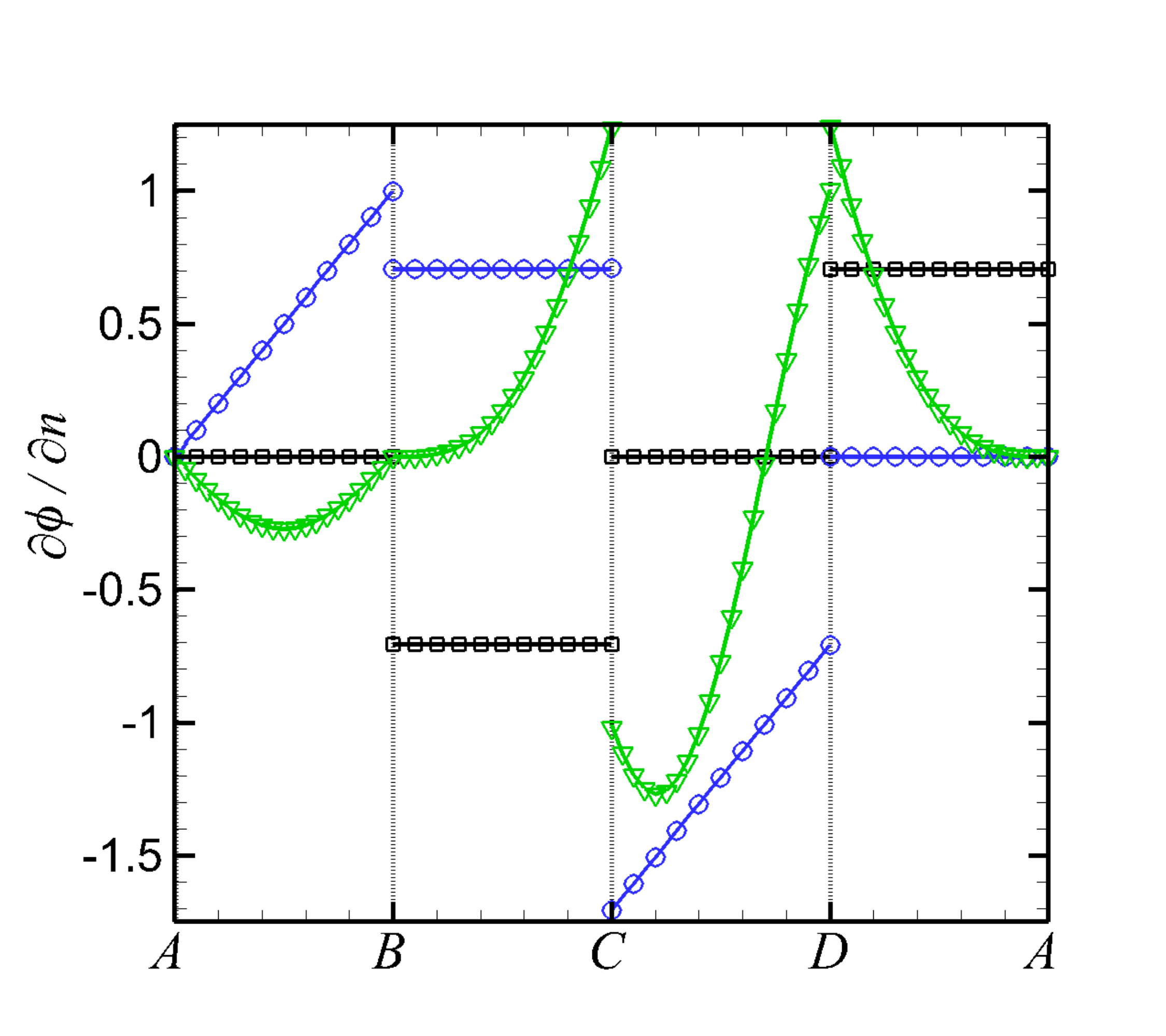} }
\subfloat[$\beta = 90^{\text{o}}$]{ \includegraphics[width=0.5\textwidth] {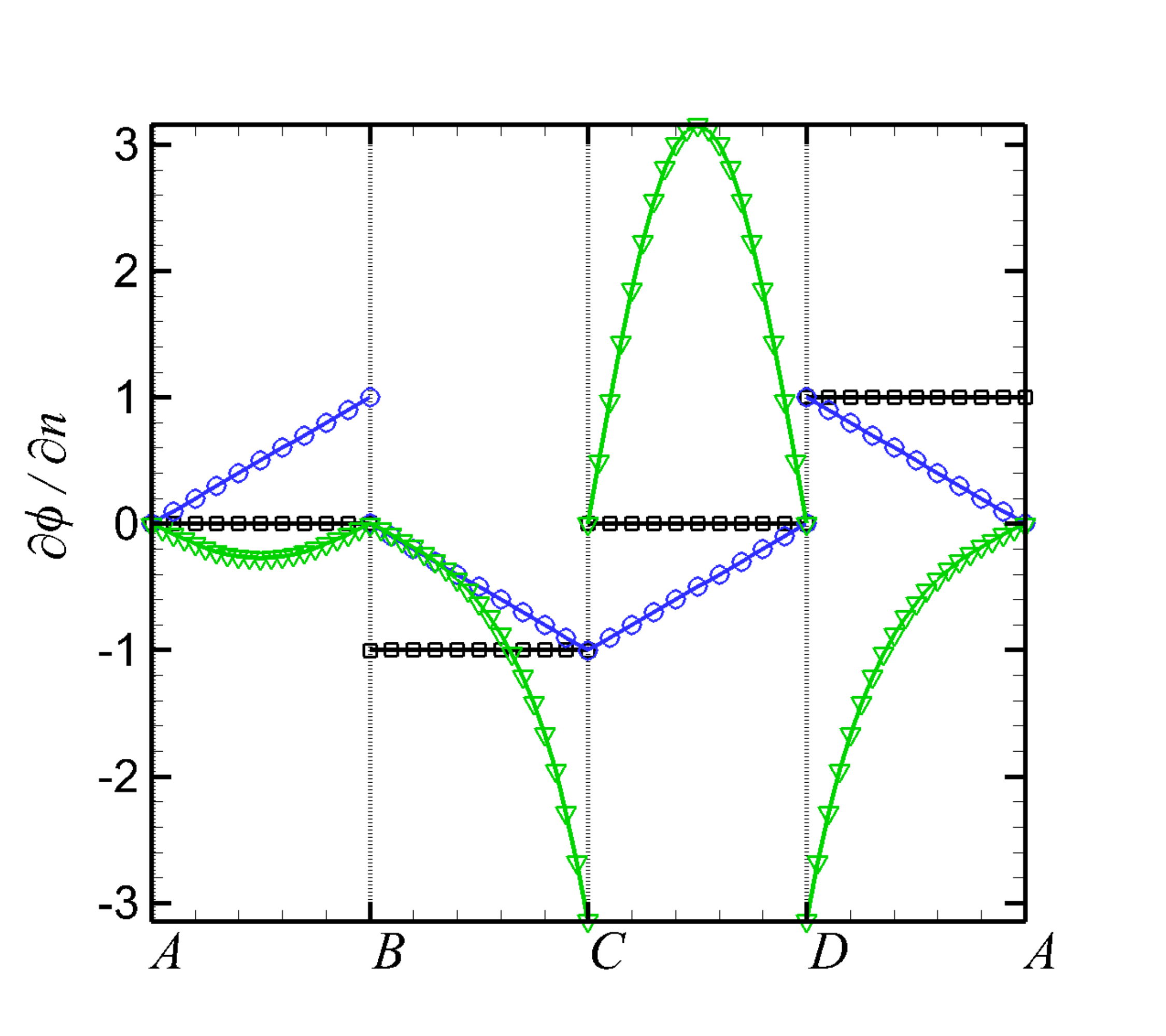} }
\caption{Variation in the normal derivative along two different domain boundaries corresponding to Fig.~\ref{fig:cornersketch} as $\beta = 45^{\text{o}}$ and $\beta = 90^{\text{o}}$, between the exact solutions given by Eqs. (19)-(21) (lines) and the results obtained by our non-singular BIM (symbols). For clarity, only a subset of the results on the 21 nodes on each side of the parallelogram are shown for Case I: black squares, Case II: blue circles and Case III: green triangles.} \label{fig:crnerr}
\end{figure}

We solved the above corner problems with Dirichlet boundary conditions. The comparisons between the results for the normal derivatives obtained by our non-singular BIM with 21 nodes on each edge and the analytical solutions are shown in Fig.~\ref{fig:crnerr} when $\beta = 45^{\text{o}}$ and $\beta = 90^{\text{o}}$. The absolute errors for these cases can be found in Table~\ref{Tbl}.

% -- TABLE 1 --
\begin{table*}\centering
    \caption{Absolute error of the 2D non-singular boundary integral mehtod for solving Cases I, II and III given by Eqs. (19)-(21) for the domain in Fig.~\ref{fig:cornersketch} for different corner angle $\beta$. }\label{Tbl}
    %\ra{1.3}
    \begin{tabular}{l c c c}
        \toprule
        Domain angle & Case I & Case II & Case III \\
        \midrule\midrule
        $\beta = 45^{\text{o}}$ & 0.0046\% & 0.29\% & 0.87\% \\
        \midrule
        $\beta = 90^{\text{o}}$ &  0.012\% & 0.019\% & 0.77\% \\
        \bottomrule
    \end{tabular}
\end{table*}

Other pairs of $f^{L}(\boldsymbol{x})$  and $f^{R}(\boldsymbol{x})$ that satisfy Eqs.~(\ref{eq:crf(x)L}) and~(\ref{eq:crf(x)R}) can also be used, for example
\begin{align}\label{eq:crnrchc2L}
f^{L}(\boldsymbol{x}) = \frac{|\boldsymbol{x}_{0}-\boldsymbol{x}^{L}_D|^2} {\boldsymbol{n}^{L}_0\bcdot(\boldsymbol{x}_{0}-\boldsymbol{x}^{L}_D)} \ln{\left(\frac{|\boldsymbol{x}-\boldsymbol{x}^{L}_D|}{|\boldsymbol{x}_{0}-\boldsymbol{x}^{L}_D|}\right)}
\end{align}
\begin{align}\label{eq:crnrchc2R}
f^{R}(\boldsymbol{x}) = \frac{|\boldsymbol{x}_{0}-\boldsymbol{x}^{R}_D|^2} {\boldsymbol{n}^{R}_0\bcdot(\boldsymbol{x}_{0}-\boldsymbol{x}^{R}_D)} \ln{\left(\frac{|\boldsymbol{x}-\boldsymbol{x}^{R}_D|}{|\boldsymbol{x}_{0}-\boldsymbol{x}^{R}_D|}\right)}
\end{align}
where $\boldsymbol{x}^{L}_D$ and $\boldsymbol{x}^{R}_D$ are located outside of the calculation domain, and
\begin{align}
\left\{
\begin{aligned}
\boldsymbol{n}^{L}_0\bcdot(\boldsymbol{x}_{0}-\boldsymbol{x}^{L}_D)\neq0, \quad \boldsymbol{n}^{R}_0\bcdot(\boldsymbol{x}_{0}-\boldsymbol{x}^{L}_D)=0; \\
\boldsymbol{n}^{R}_0\bcdot(\boldsymbol{x}_{0}-\boldsymbol{x}^{R}_D)\neq0, \quad \boldsymbol{n}^{L}_0\bcdot(\boldsymbol{x}_{0}-\boldsymbol{x}^{R}_D)=0.
\end{aligned}
\right.
\end{align}
This pair of $f^{L}(\boldsymbol{x})$ and $f^{R}(\boldsymbol{x})$ are also used to solve the above three cases with Dirichlet boundary conditions when $\beta = 90^{\text{o}}$. Once again, 21 nodes are employed on each edge. The absolute errors between the results for the normal derivatives obtained by our non-singular BIM and the analytical solutions for these cases are the same as Table~\ref{Tbl}.

This method of treating corner nodes can be extended to handle edges and vertices in 3D domains, but will not be considered here.

\subsection{Wave drag at a semi-infinite deformable boundary}
\begin{figure}[t]
\centering
\includegraphics[width=0.7\textwidth] {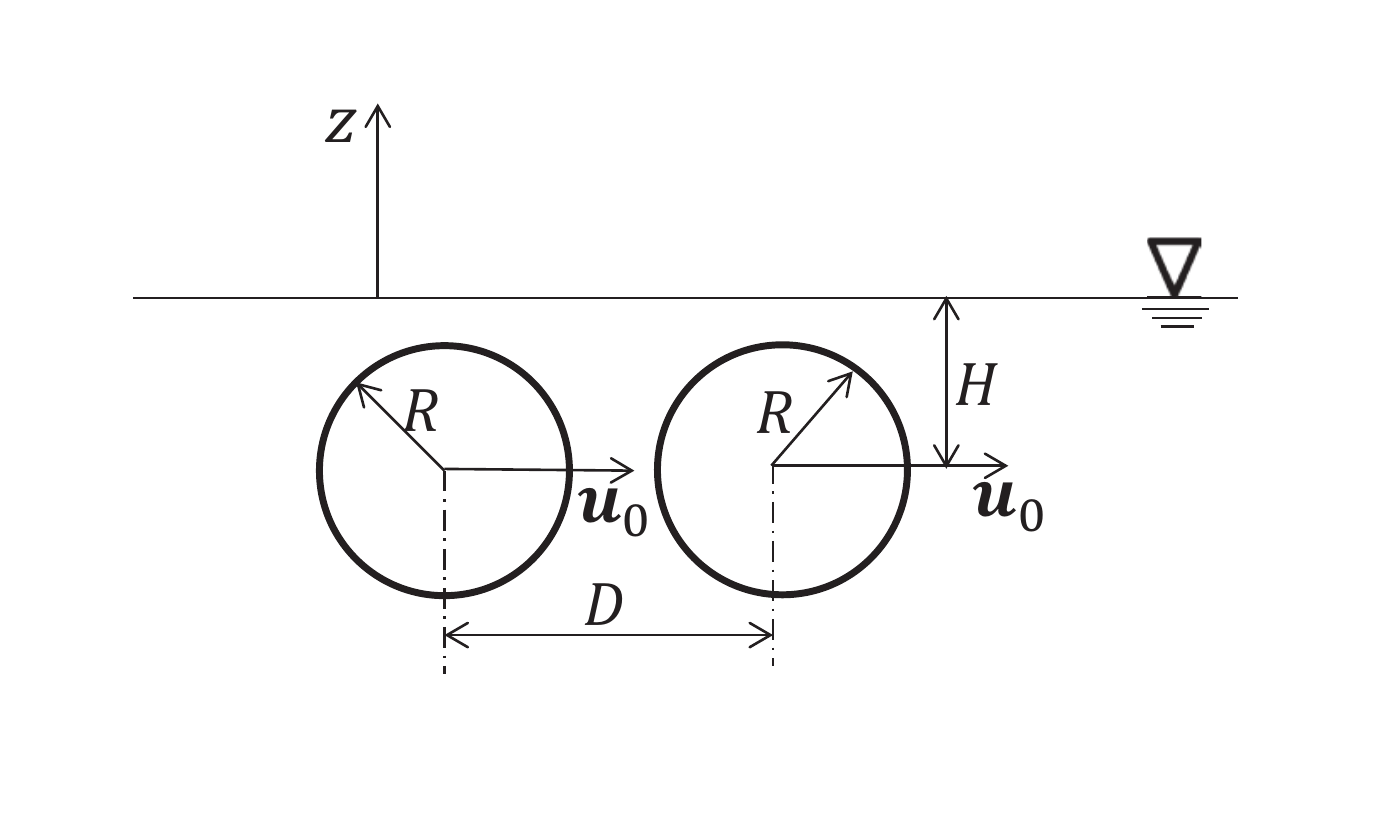}
\caption{The geometry of two spheres translating near a deformable free interface.}\label{fig:illus_internalxd}
\end{figure}

To illustrate the generality of our non-singular formulation of the BEM, we examine the following wave drag problem that has sufficient complexity to be interesting~\citep{Landau1966}. Consider two identical spheres of radius, $R$, separated by a constant distance $D$ between their centres and moving with constant speed $\boldsymbol{u}_0 = U_0\boldsymbol{i}$ at the same depth $H$ below an infinite deformable free surface in a gravity field, as shown in Fig.~\ref{fig:illus_internalxd}. The origin of the global reference frame is set to coincide with the initially undisturbed free surface (at $z=0$) with the surface elevation $z$ pointing upwards (Fig.~\ref{fig:illus_internalxd}). For simplicity, the upper phase (air) is assumed to have negligible mass density and interfacial tension effects have been omitted, although it is easy to dispense with such simplifications at the expense of introducing more physical parameters. This is a time-dependent potential flow problem as surface waves will be generated on the deformable free surface while the pair of spheres travelling beneath it.

% -- FIGURE 5 --

The potential flow velocity field, $\boldsymbol{u} = \bnabla \phi$, generated by the moving spheres is found by solving the Laplace equation for the velocity potential, $\phi$ at each time step. On the spheres, the velocity is given and the potential is calculated. Initially, $\phi = 0$ on the undisturbed free surface. Solving the mixed boundary value problem defined by the free surface and the two spheres, provides values of the normal velocity $\partial{\phi}/\partial{n}$ at the free surface that is then used to predict its shape at the next time step. The value of $\phi$ on the free surface at the next time step is found from the unsteady Bernoulli equation evaluated on the surface: $\rho \text{D}\phi/\text{D}t = \rho u^2 /2 - \rho g  z$ where $\text{D}/\text{D}t = \partial{}/\partial{t} + \boldsymbol{u}\cdot\bnabla$ is the material derivative, $\rho$ the fluid density and $g$ the gravitational acceleration. Such time-stepping then gives the spatio-temporal evolution of the surface waves. The wave drag force experienced by the two spheres that are moved at constant velocity in this 3D problem is calculated by integrating the pressure $p$ on the sphere surfaces.

The deformable interface extends to infinity in the $xy$-direction where it asymptotes to a flat surface. Therefore far from the spheres, the potential and its derivative vanish asymptotically like $1/|\boldsymbol{x}|$ and $1/|\boldsymbol{x}|^2$ as $|\boldsymbol{x}|\rightarrow\infty$. We therefore assume that beyond a radius of 40$R$ from the spheres, we can take the interface to be flat and the potential and its derivative also take on their asymptotic forms so that we can evaluate analytically the contribution to the surface integrals from the far field, see for example~\citep{Gao1998, Ribeiro2009}. Since Eq.~(\ref{eq:nsbim_f}) is non-singular, there is no need to map the far field elements to two triangular elements or introduce any artificial points at infinity as was done in~\cite{Ribeiro2009}. Finally, the integral over the half spherical surface at infinity will contribute two terms $2\pi \phi(\boldsymbol{x}_0)$ and $2\pi (\partial{\phi}/\partial{n})_0  |\boldsymbol{x}_0 - \boldsymbol{x}_D|^2  / [\boldsymbol{n}_0 \bcdot (\boldsymbol{x}_0 - \boldsymbol{x}_D)]$ to the left-hand side of Eq.~(\ref{eq:nsbim_f}).

% -- FIGURE 6 --
\begin{figure}[t]
\centering
\subfloat[]{ \includegraphics[width=0.5\textwidth] {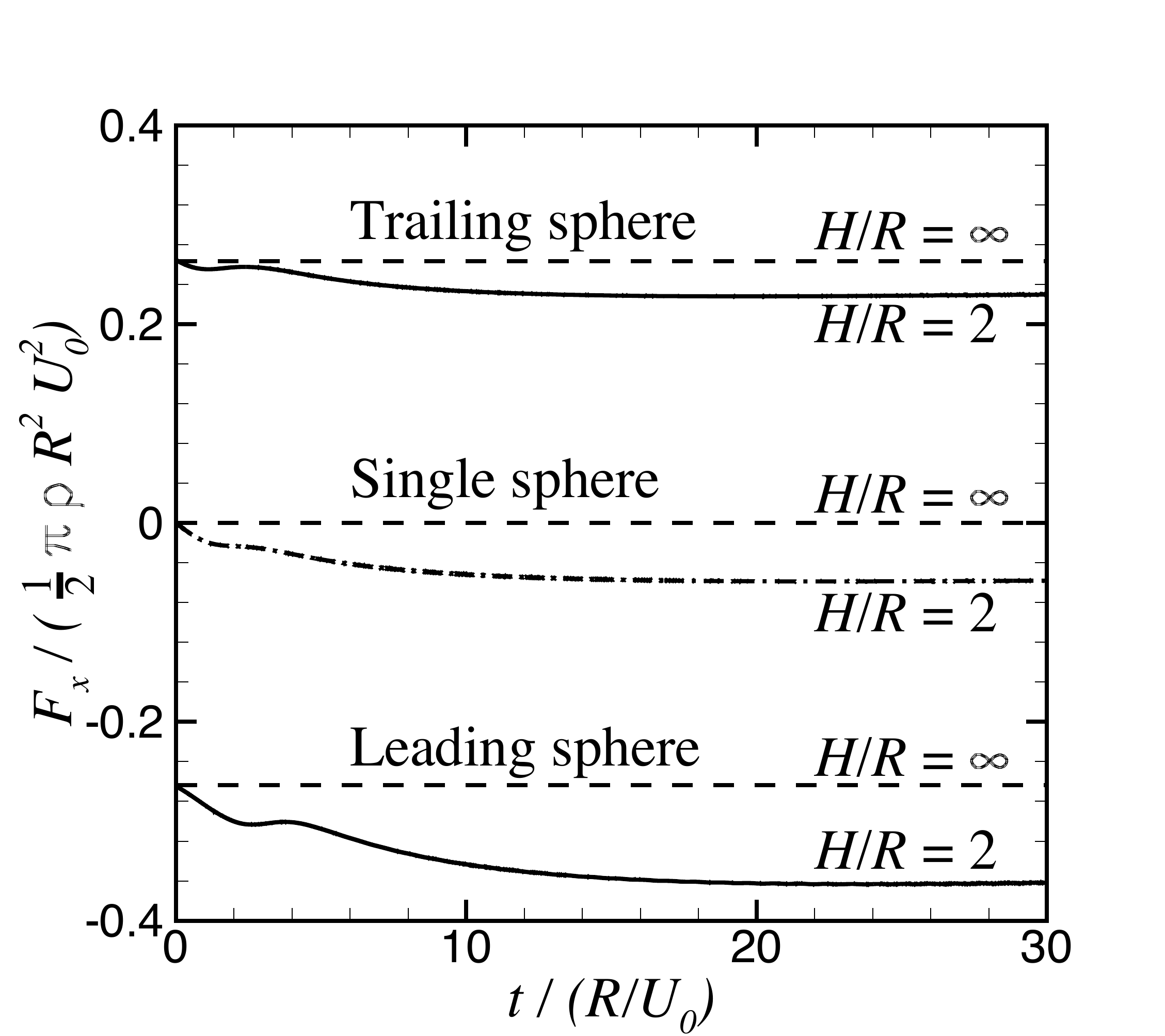} }
\subfloat[]{ \includegraphics[width=0.5\textwidth] {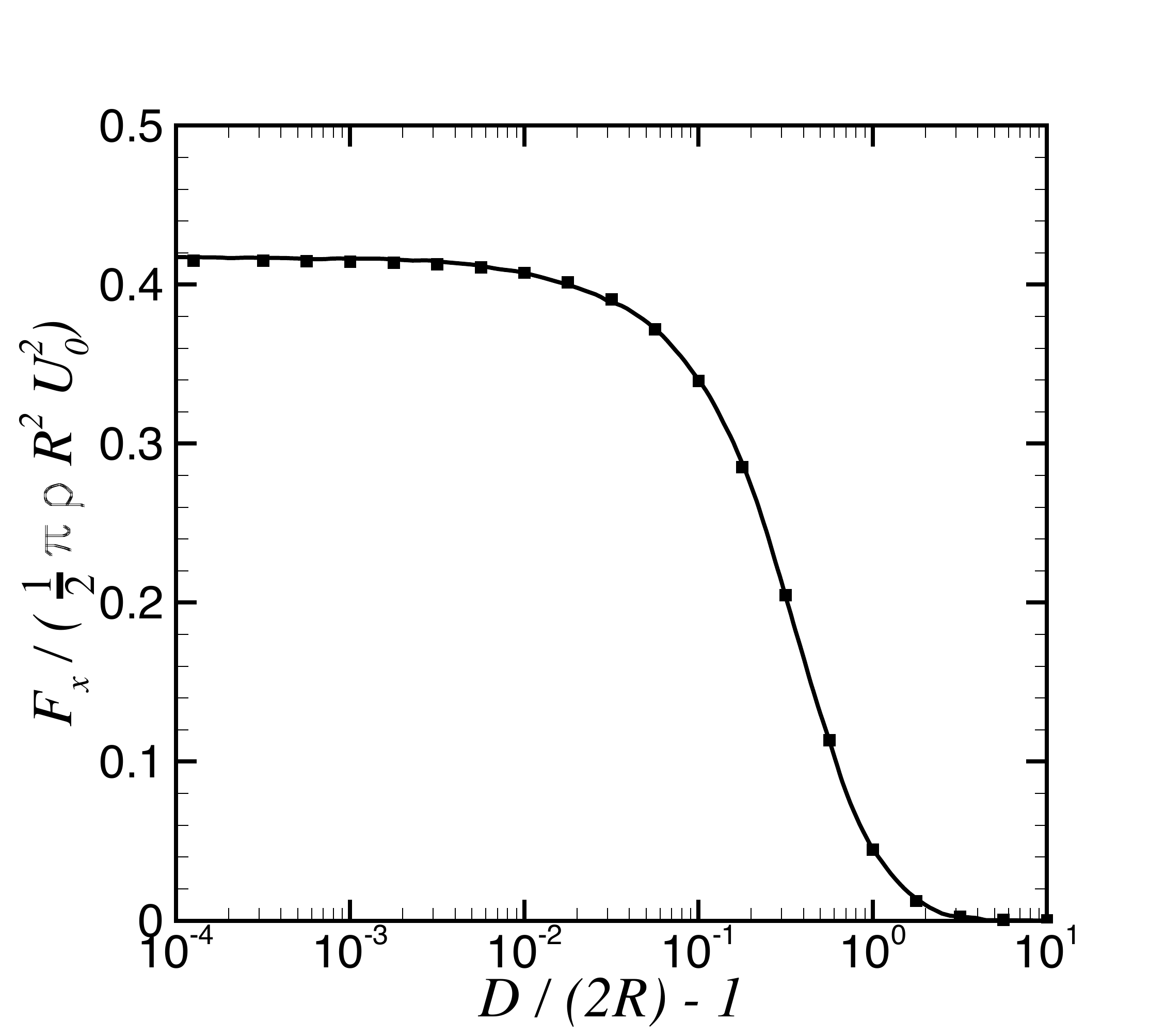} }
\caption{(a) The time dependent force, $\boldsymbol{F}=F_x \boldsymbol{i}$ experienced by a pair spheres of radius $R$ at separation, $D = 2.4 R$, moving at constant velocity, $\boldsymbol{u}_0=U_0 \boldsymbol{i}$ from rest along the line of centres at different depths, $H$, parallel to a deformable free surface in a gravity field at Froude number, $Fr = U_0/\sqrt{gH} = 1$. Corresponding results for a single sphere are given for comparison. (b) Variations of the force on the trailing sphere with separation, $D$ in the absence of the interface ($H/R = \infty$) obtained from a 3D calculation using the present non-singular BEM (points) and from the analytic result of \cite{Miloh1977} (line). The forces are scaled to give the usual drag coefficient, $C_d$. } \label{fig:forces}
\end{figure}

On each sphere, 1280 linear triangular elements with 642 nodes are employed. On the free surface, 15000 elements with 7651 nodes are used. When $\boldsymbol{x}_0$ is on a sphere, $\boldsymbol{x}_D$ in Eq.~(\ref{eq:psi_inverse}) is taken to be the centre of the sphere. When $\boldsymbol{x}_0$ is on the free surface, we choose the Cartesian components of $\boldsymbol{x}_D = (x_D, y_D, z_D)$ to be $z_D = 3 R$ and any convenient $x_D$ and $y_D$ that ensure $\boldsymbol{n}_0 \bcdot (\boldsymbol{x}_0 - \boldsymbol{x}_D) \neq 0$.

With a constant time step of  $0.02R/U_0$, the position of each surface node is updated with $\text{D}\boldsymbol{x}/\text{D}t = \boldsymbol{u} \equiv \bnabla \phi (\boldsymbol{x})$ and the spheres are translated using $\boldsymbol{u}_{0}$. The calculation continues until the spheres have travelled a distance of $30R$. The elastic mesh technique~\citep{Wang2003} is applied on the free surface to ensure a uniform mesh even after many time-steps.

The drag force acting on each sphere is calculated by $\boldsymbol{F}=\int_{S}p\text{ }\boldsymbol{n}\text{ d}S$ where $S$ is the surface of the sphere. In Fig.~\ref{fig:forces}a we present the time variation of the force $\boldsymbol{F} = F_x \boldsymbol{i}$ on each of the two spheres travelling with constant velocity $\boldsymbol{u}_0 = U_0 \boldsymbol{i}$. The three non-dimensional parameters that govern this problem are: $H/R$, $D/R$ and the Froude number, $Fr = U_0/\sqrt{gH}$ (Fig.~\ref{fig:illus_internalxd}). For our illustrative example, we chose a sphere separation of $D/R = 2.4$, a submerged depth of $H/R = 2$,  and $Fr = 1$. Also shown are results for $H/R = \infty$ that corresponds to the absence of the deformable interface and the force on an isolated single sphere that corresponds to  $D/R = \infty$. The leading sphere experiences a retarding force in the $x$-direction whereas the trailing sphere experiences a force in the direction of travel.  At $H/R = \infty$, these two forces are equal and opposite as expected from the d'Alembert Paradox of potential flow. The proximity of the moving spheres to the deformable interface provides a net non-zero wave drag on the pair of moving spheres~\citep{Landau1966}. We point out that the treatment of the initial condition of this problem has been simplified by omitting the detailed effects of acceleration from rest that will give rise to an added mass term. The inclusion of this effect complicates the equation of motion, but does not affect the demonstration of the utility of the present non-singular formulation of the BIM.

In Fig.~\ref{fig:forces}b we show the variation of the drag force on the trailing sphere as a function of the sphere separation, $D$ at $H/D = \infty$. There is excellent agreement between our non-singular boundary integral method and the analytic result of~\cite{Miloh1977}.

\begin{figure}[!t]
\centering
\subfloat[]{ \includegraphics[width=0.47\textwidth] {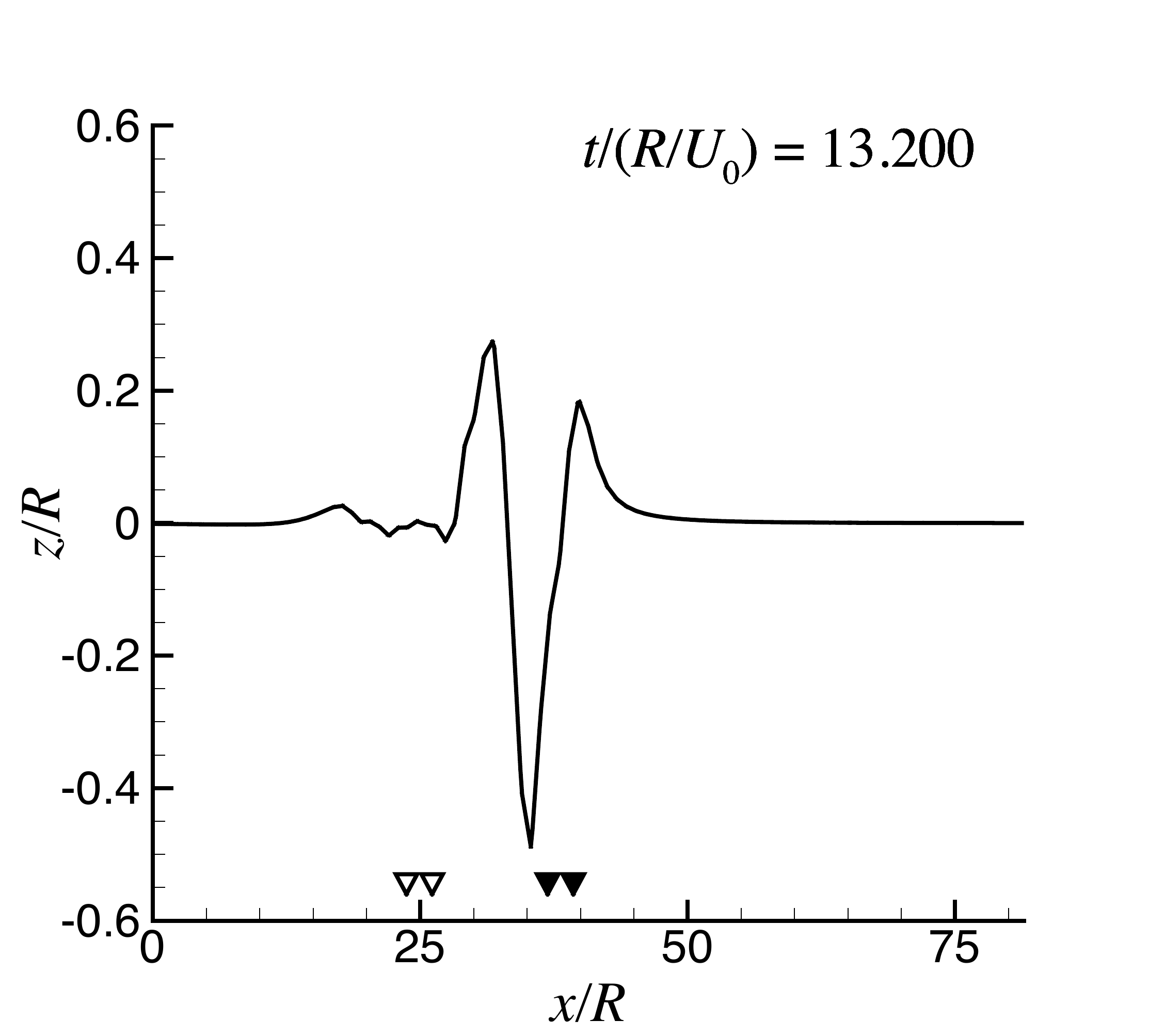} }
\subfloat[]{ \includegraphics[width=0.47\textwidth] {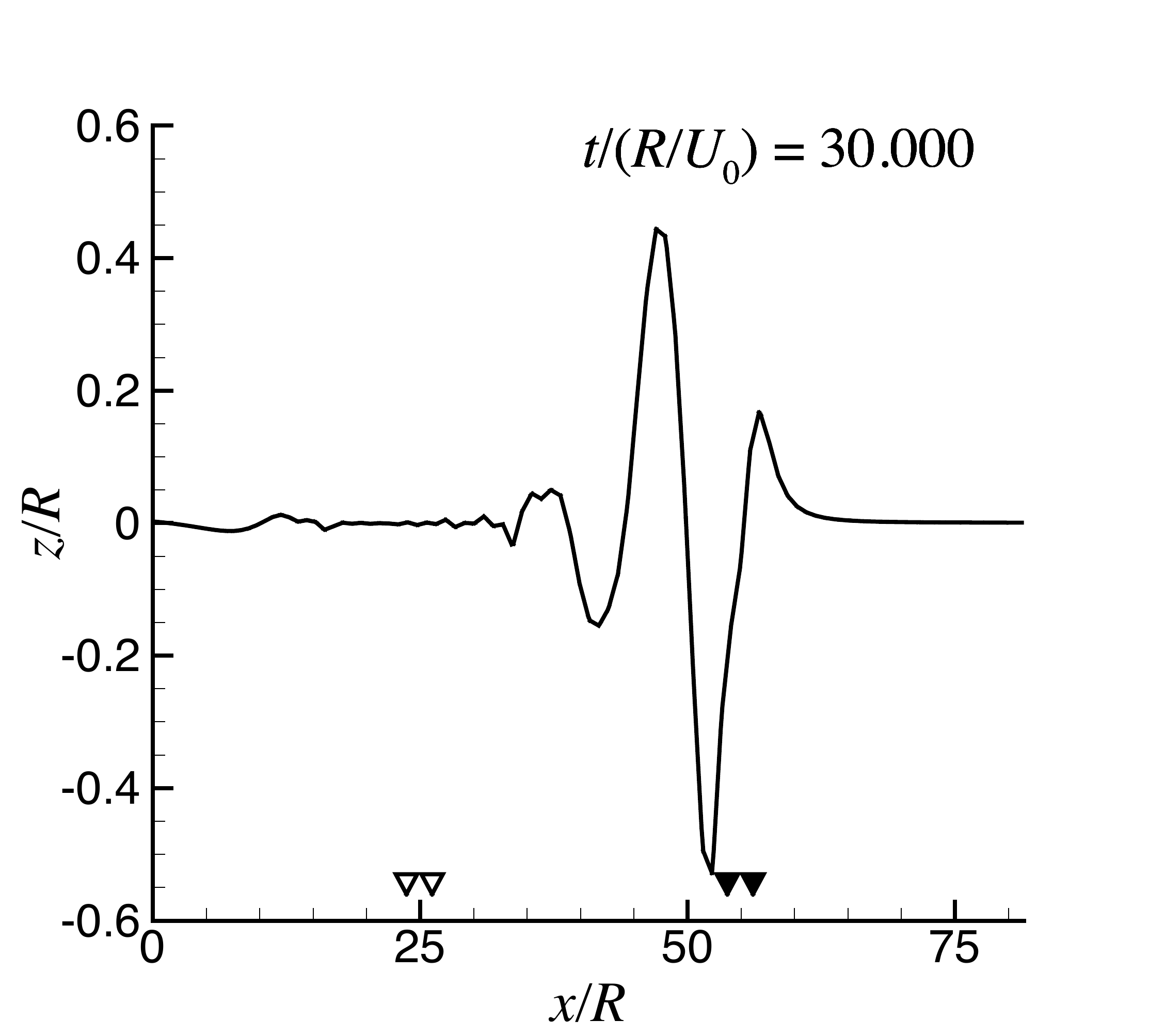} }\\
\subfloat[]{ \includegraphics[width=0.43\textwidth] {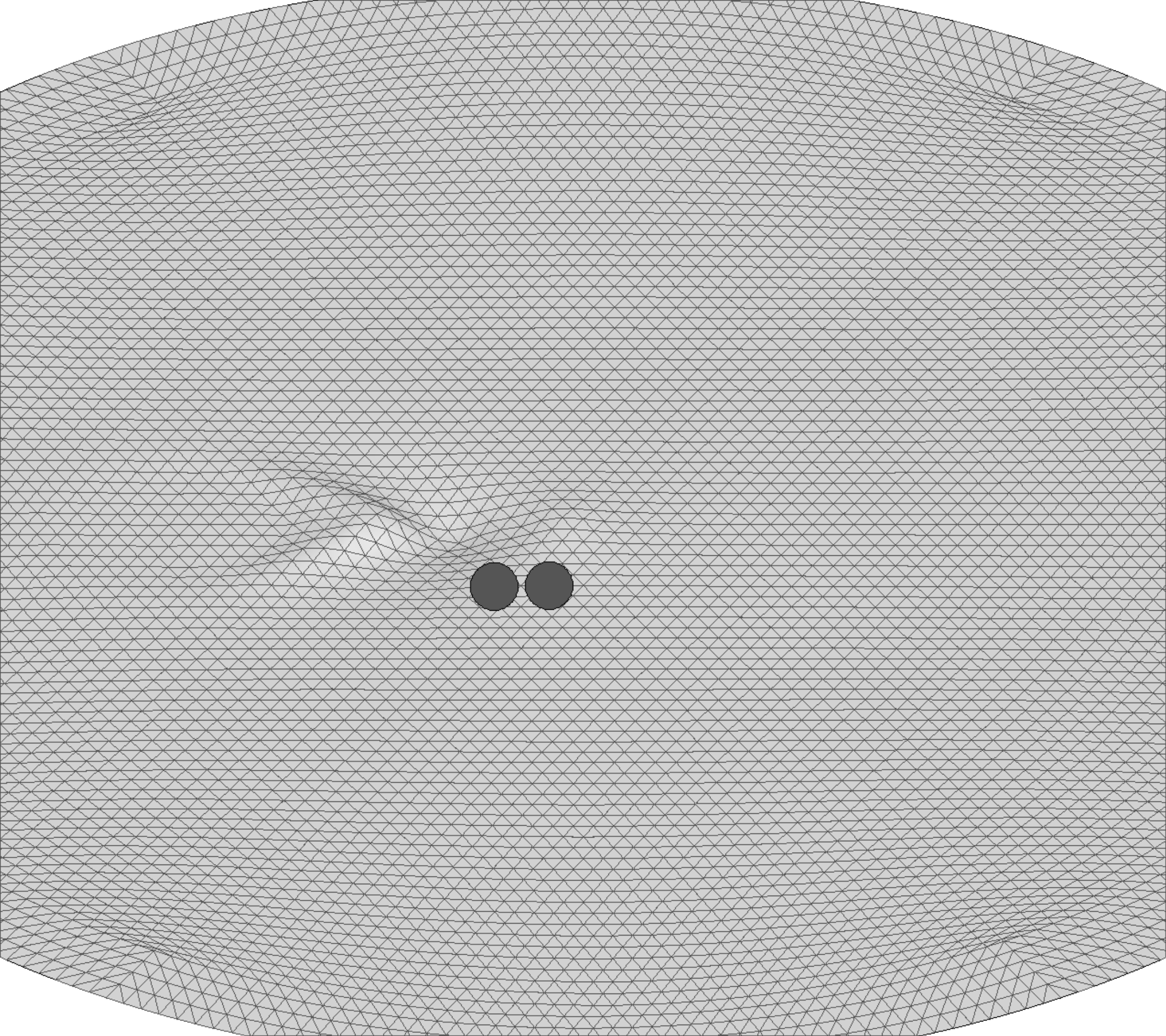} }
\subfloat[]{ \includegraphics[width=0.43\textwidth] {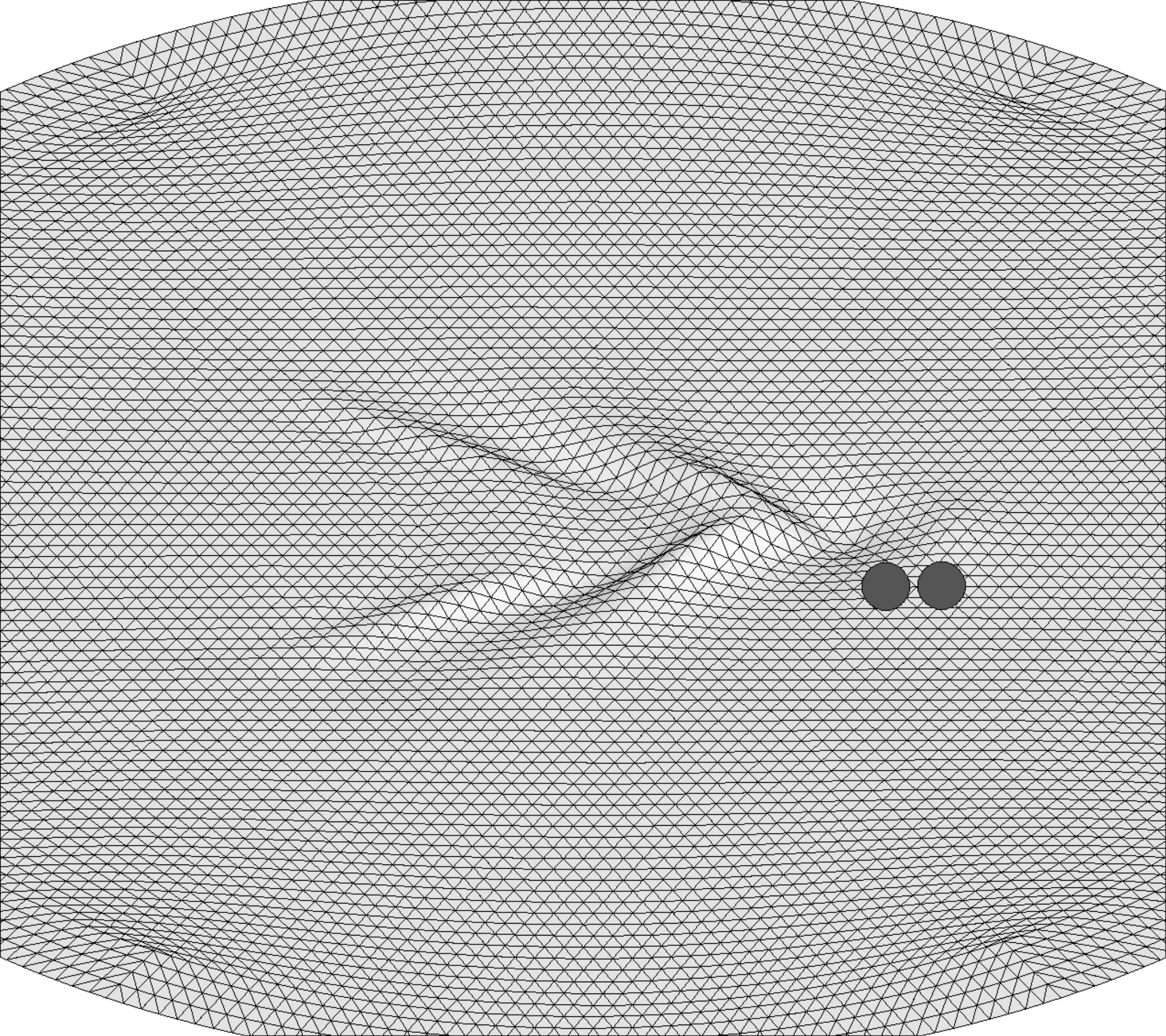} }
\caption{Variation of the surface wave amplitude generated by two submerged spheres initially at positions indicated by the open triangles, as they travel along the $x$-direction. The positions of the spheres at times: (a) $t = 13.20 R/U_0$ and (b) $t = 30.00 R/U_0$ are indicated by solid triangles. The corresponding 3D representations of the surface waves are shown in (c) and (d). The parameter values are the same as for Fig.~\ref{fig:forces}a. Animations corresponding to these figures are available as online e-component material.} \label{fig:frmvmnt}
\end{figure}

Snapshots of the spatio-temporal variations of the surface waves generated by the pair of moving submerged spheres are shown in Fig.~\ref{fig:frmvmnt}. Variations of the surface wave amplitude along the direction of travel at two different times are shown in Fig.~\ref{fig:frmvmnt} (a) and (b). Three dimensional representations of the corresponding surface waves in Fig.~\ref{fig:frmvmnt} (c) and (d) show the lateral extent of surface disturbance that obviously depends on the sphere spacing and the depth of immersion. Movies of the wave amplitude and the surface wave corresponding to these figures are available as online supplementary material. Although potential flow is conservative with no energy dissipation or damping mechanism, the drag waves on the free surface appear damped because of kinetic energy of the moving sphere being distributed into the infinite fluid domain~\citep{Landau1966}.

% -- FIGURE 7 --

% ---- Section 7 - Conclusions ----

\section{Conclusions}

In this paper, we have provided a fully non-singular formulation of the boundary integral method for potential problems. The usual singularities in the kernels that arise from the fundamental solutions are removed analytically without introducing additional unknowns or extra equations to be solved. Apart from excising the ``mathematical monsters'' in the conventional boundary integral formulation~\citep{Becker92}, the amount of computer code needed to implement the non-singular boundary integral equation has been reduced by about 60\%. For the special case of axisymmetric problems, our formulation also removed the technical inconvenience associated with nodes on the axis of symmetry~\citep{Becker92}. The robustness and generality of our approach has been demonstrated with examples involving osculating boundaries, domains with corners and a wave drag due to an infinite deformable boundary.

\section*{Acknowledgement}
DYCC is a Visiting Scientist at the Institute of High Performance Computing and an Adjunct Professor at the National University of Singapore. This work is supported in part by the Australian Research Council Discovery Project Grant Scheme.



%% References


\bibliographystyle{model1-num-names}
\bibliography{Pot_NSBIM}






% \end{thebibliography}

\end{document}